\documentstyle[twoside,12pt]{article}

\newcommand{\chapter}{\section}

\begin{document}

\newtheorem{Thm}{Theorem}
\newtheorem{Ax}{Axiom}
\newtheorem{Prop}{Proposition}
\newtheorem{Cor}[Prop]{Corollary}
\newtheorem{Main}{}
\renewcommand{\theMain}{}
\newtheorem{Lem}[Prop]{Lemma}
\newtheorem{Fact}{Fact}
\renewcommand{\theFact}{}

\newtheorem{Def}{Definition}
\newtheorem{rmk}{Remark}
\newenvironment{Rmk}{\begin{rmk}\em}{\end{rmk}}
\newtheorem{exm}{Example}
\newenvironment{Exm}{\begin{exm}\em}{\end{exm}}

\newcommand{\qed}{\par {\EM QED} }
\newtheorem{prf}{Proof}
\renewcommand{\theprf}{}
\newenvironment{Prf}{\begin{prf}\em}{\qed\end{prf}}
\newtheorem{prff}{}
\renewcommand{\theprff}{}
\newenvironment{Prff}{\begin{prff}\em}{\qed\end{prff}}





\newcommand{\YES}[1]{#1}
\newcommand{\NOT}[1]{}

\newcommand{\cA}{{\cal A}}
\newcommand{\cB}{{\cal B}}
\newcommand{\cC}{{\cal C}}
\newcommand{\cD}{{\cal D}}
\newcommand{\cE}{{\cal E}}
\newcommand{\cF}{{\cal F}}
\newcommand{\cG}{{\cal G}}
\newcommand{\cH}{{\cal H}}
\newcommand{\cI}{{\cal I}}
\newcommand{\cJ}{{\cal J}}
\newcommand{\cK}{{\cal K}}
\newcommand{\cL}{{\cal L}}
\newcommand{\cM}{{\cal M}}
\newcommand{\cN}{{\cal N}}
\newcommand{\cO}{{\cal O}}
\newcommand{\cP}{{\cal P}}
\newcommand{\cQ}{{\cal Q}}
\newcommand{\cR}{{\cal R}}
\newcommand{\cS}{{\cal S}}
\newcommand{\cT}{{\cal T}}
\newcommand{\cU}{{\cal U}}
\newcommand{\cV}{{\cal V}}
\newcommand{\cW}{{\cal W}}
\newcommand{\cX}{{\cal X}}
\newcommand{\cY}{{\cal Y}}
\newcommand{\cZ}{{\cal Z}}

\newcommand{\bbb}[1]{{\mbox{\bf #1}}}

\newcommand{\bN}{\bbb{N}}
\newcommand{\bZ}{\bbb{Z}}
\newcommand{\bR}{\bbb{R}}
\newcommand{\bC}{\bbb{C}}
\newcommand{\bQ}{\bbb{Q}}
\newcommand{\bT}{\bbb{T}}

\newcommand{\noind}[1]{{\setlength{\parindent}{0cm} #1}}
\newcommand{\parsk}{\par\medskip}

\newcommand{\varend}{\end{document}}

\newcommand{\LP}{\left(}  \newcommand{\RP}{\right)}
\newcommand{\LQ}{\left[}  \newcommand{\RQ}{\right]}
\newcommand{\LB}{\left\{} \newcommand{\RB}{\right\}}
\newcommand{\LA}{\left\langle} \newcommand{\RA}{\right\rangle}
\newcommand{\LS}{\left|}  \newcommand{\RS}{\right|}
\newcommand{\LN}{\left\|} \newcommand{\RN}{\right\|}
\newcommand{\bLP}{\bigl(}  \newcommand{\bRP}{\bigr)}

\newcommand{\ts}{{\otimes}} \newcommand{\TS}{{\bigotimes}}

\newcommand{\es}{\emptyset}
\newcommand{\sm}{\setminus} \newcommand{\st}{\subset}
\newcommand{\eps}{\varepsilon}

\newcommand{\beware}{~$\Psi\!\Psi\!\Psi$~}

\newcommand{\EM}{\em\bf}

\newcommand{\BE}{\begin{equation}}  \newcommand{\EE}{\end{equation}}
\newcommand{\BER}{$$\begin{array}}  \newcommand{\EER}{\end{array}$$}
\newcommand{\head}[1]{
             \par\bigskip\noind{{\large\bf#1}\nopagebreak\par}\medskip}
\newcommand{\chead}[1]{\begin{center}
              \par\bigskip{\large\bf#1}\nopagebreak\par\medskip
              \end{center}}

\newcommand{\fillitem}{\L{\embox{}}\vspace{-.6cm}}

\newcommand{\dfrac}[2]{\displaystyle{\frac{#1}{#2}}}
\newcommand{\toprove}[1]{{\buildrel\mbox{?}\over#1}}
\newcommand{\DEF}{{\buildrel\mbox{\scriptsize\,\,def\,\,}\over=}}
\newcommand{\BAR}[1]{\overline{#1}}
\newcommand{\I}{\infty}
\newcommand{\al}{\alpha} \newcommand{\be}{\beta}
\newcommand{\OM}{\Omega} \newcommand{\om}{\omega}
\newcommand{\la}{\lambda}
\newcommand{\mod}{\mbox{ mod }}
\newcommand{\half}{\frac{1}{2}}
\newcommand{\NI}{{n\to\infty}}
\newcommand{\ang}{{<\!\!\!)}}
\newcommand{\arc}[1]{\breve{#1}}

\newenvironment{txeqn}[1]{\begin{equation}\label{#1}\begin{array}{l}}{
                         \end{array}\end{equation}}



\newcommand{\cs}{continuous}
\newcommand{\fn}{function}
\newcommand{\bd}{bounded}
\newcommand{\nd}{neighborhood}
\newcommand{\tg}{topolog}

\title{Topological Interpretation of Function Spaces Stable under a General
Operation}
\author{Eliahu Levy\\
Department of Mathematics\\
Technion -- Israel Institute of Technology,
Haifa 32000, Israel\\
email: eliahu@techunix.technion.ac.il}


\date{}


\maketitle
\begin{abstract}
Function (linear) spaces $E$ on which some arbitrary real or complex
\cs\ \fn\ operates (in other words, $E$ is stable under the unary operation
defined by this \fn) were studied by K.~de Leeuw, Y.~Katznelson, Y.~Sternfeld
and Y.~Weit, showing, e.g.\ that for a complex \fn\ space $E$ consisting of
bounded \fn s, containing the constants and closed
w.r.t.\ the $\sup$-norm, if some non-affine \cs\ \fn\ operates on $E$
then every analytic \fn\ operates on $E$: $E$ is an algebra, and if further
some non-analytic \cs\ \fn\ (say, complex conjugation) operates on $E$ then any
\cs\ \fn\ does so: $E$ is a self-adjoint algebra.

In this note we approach this issue from a different point of view.
Our scalars are any field $K$, and the \fn s are defined on an abstract set
$I$, with both $K$ and $I$, at the outset, without a \tg y.
It turns out that the very fact that the \fn\ space $E$ is stable w.r.t.\ a
non-affine operation (no continuity required) allows it to induce a topology
on $I$, while the space of \fn s which operate on it induces \tg ies
on $K^n$, $n=1,2,\ldots$. The \tg ized $I$ is used to prove a density
property for such $E$. Also, these topologies allow us, in some cases, to
investigate general ``homomorphisms'' between \fn\ spaces w.r.t.\ to
non-affine operations under which they are stable.
\end{abstract}

\setcounter{section}{-1}

\section{Definitions}

Let $K$ be a field\NOT{(everything works also for a division ring?)}.
We refer to any $f:K^n\to K$ as an {\EM operation} and say that it is
{\EM additively affine} if it can be written as a sum of a constant and a homomorphism
of the additive groups.
\parsk

Let $I$ be a set, and let $E\subset K^I$ be a (linear) subspace (a \fn\ space).
We call $E$ {\EM unital} if it {\em contain the constant \fn s}, which we
always assume. Let $f:K^n\to K$ be an operation.
We say that {\EM $E$ is stable w.r.t.\ the operation $f$}
or that $f$ {\EM operates on $E$},
if for any $\xi_1,\ldots,\xi_n\in E$ also
$t\in I\mapsto f\LP\xi_1(t),\ldots,\xi_n(t)\RP$ is in $E$.
We say that $E$ is stable w.r.t.\ a set of \fn s (operations) if it is
stable w.r.t.\ each member.
\parsk

We shall say that $E$ {\EM separates points} in $I$ if for any different
$t_1,t_2\in I$ and any $a_1,a_2\in K$ $\exists$ a $\xi\in E$
with $\xi(t_1)=a_1$,\,$\xi(t_2)=a_2$. Equivalently, for no different
$t_1$ and $t_2$, $\xi(t_1)$ and $\xi(t_2)$ are proportional for the $\xi\in E$.

\section[Induced Topology]
{Function spaces stable under a non-additively-affine operation induce a topology}

We try to make $E$ define a \tg y on $I$. Call a subset $S\subset I$
{\EM closed} if for any $t\in I\setminus S$ $\exists$ a $\xi\in E$ vanishing
on $S$ but not at $t$. It is clear that $\es$ and $I$ are closed and that
any intersection of closed sets is closed (if $t\notin\cap_\al S_\al$ then
$t\notin S_\al$ for some $\al$ and a $\xi$ as above for $S_\al$ will serve
also for $\cap_\al S_\al$). Closed sets can also be defined as sets of
solutions of (possibly infinite) sets of equations of the form $\xi(t)=0$
with $\xi\in E$.

\begin{Def}
Suppose the subspace $E\subset K^I$ is {\em unital}, i.e.\ contains the
constants. We say that $E$ is {\EM \tg y-inducing} if the union of any two
closed subsets of $I$ (defined as above) is closed, so one has a \tg y on
$I$, which we call {\EM the \tg y induced by $E$}.
\end{Def}
\parsk

Note that if $E$ is \tg y-inducing and separates points, then the
induced \tg y is $T_1$ (i.e.\ every singleton is closed).
\parsk

Note also that a general subspace need not be \tg y-inducing. For example,
for the hyperplane 
$$E=\LB\xi\in K^{\{1,\ldots,n\}}\,\Big|\,\sum_{i=1}^n\lambda_i\xi_i=0\RB$$
with all $\lambda_i\ne0$ and $\sum_i\lambda_i=0$ (making $E$ unital), the
closed sets are the subsets $S\subset\{1,\ldots,n\}$
whose number of elements is different from $n-1$.
\parsk

Of course, for $I=K^m$, $E=\{$the polynomials$\}$ is \tg y-inducing
and we obtain the {\em Zariski \tg y}, while for $K=\bR$ (the reals), $I$
a Hausdorff compact \tg ical space and $E=\{$the \cs\ \fn s$\}$ we obtain
the original \tg y on $I$.
\parsk

The following assertion is immediate
\begin{Prop} \label{I'}
If $E$ is \tg y-inducing on $I$ and $I'\subset I$, then the space of all
the restrictions of members of $E$ to $I'$ is \tg y-inducing and induces
on $I'$ the relative \tg y from the \tg y that $E$ induces on $I$.
\end{Prop}\qed

\begin{Rmk}
Suppose $E$ \tg y-inducing. Note that if one takes in $K$ the $T_1$
\tg y where the closed sets are just the finite ones and the whole space,
then all members of $E$ will be \cs, and for this \tg y in $K$\,\,$I$ is
(mock!) ``completely regular'' -- for any closed set and a point not in it
there is a \cs\ \fn\ (even a member of $E$) vanishing on the closed set but
not at the point.
\end{Rmk}

\begin{Thm} \label{U}
Suppose the subspace $E\subset K^I$ contains the constants (i.e.\ is unital)
and is stable w.r.t.\ some $f:K^n\to K$ which is {\em not additively affine}. Then $E$
is \tg y-inducing on $I$.
\end{Thm}

\begin{Prf}
Let $S,T\subset I$ be closed and we wish to prove $S\cup T$ closed. 
Let $t_0\in I\setminus(S\cup T)$. We have $\xi,\eta\in E$ such that
$\xi$ vanishes on $S$ and $\eta$ on $T$ while both are different from $0$
at $t_0$. We have to show that $\exists$ a \fn\ in $E$ vanishing on
$S\cup T$ but not at $t_0$.
\parsk

Let $a=(a_1,\ldots a_n),\;b=(b_1,\ldots b_n)\in K^n$. Then
\begin{eqnarray*}
&&f(a\xi):=\LQ t\mapsto f\LP a_1\xi(t),\ldots,a_n\xi(t)\RP\RQ\in E\\
&&f(b\eta):=\LQ t\mapsto f\LP b_1\eta(t),\ldots,b_n\eta(t)\RP\RQ\in E
\end{eqnarray*}
and similarly
$$f(a\xi+b\eta)\in E.$$
From $\xi|_S=0$ and $\eta|_T=0$ it follows immediately that
\BE f(a\xi+b\eta)=f(a\xi)+f(b\eta)-f(0){\bf1}\quad\mbox{ on }{S\cup T}
\label{(fST)}\EE
where $\bf1$ is the constant \fn\ $1$.
\parsk
The difference of the two sides of $(\ref{(fST)})$ vanishes on $S\cup T$.
If it does not vanish at $t_0$, we are done. If it does vanish we have:
$$f(a\xi+b\eta)(t_0)=f(a\xi)(t_0)+f(b\eta)(t_0)-f(0)$$
that is
$$f\LP a\xi(t_0)+b\eta(t_0)\RP=f\LP a\xi(t_0)\RP+f\LP b\eta(t_0)\RP-f(0)$$
Note that $\xi(t_0)\ne0$,\,\,$\eta(t_0)\ne0$. Now let $x,y\in E^n$.
Taking $a=\LQ\xi(t_0)\RQ^{-1}x$, $b=\LQ\eta(t_0)\RQ^{-1}y$, one gets:
$$f(x+y)=f(x)+f(y)-f(0)\qquad x,y\in E^n,$$
making $f$ additively affine, a contradiction.
\end{Prf}

\section{A density property}

For \tg y-inducing \fn\ spaces, in particular for those unital and stable
for a non-additively-affine operation, we have the following density property:

\begin{Thm}
Let $E\subset K^I$ be separating points and \tg y-inducing. If $I$
is finite Then $E$ is the whole $K^I$. Consequently for general $I$ for any
finite subset $I'\subset I$ all \fn s on $I'$ are obtained as restrictions of
members of $E$ (density property).
\end{Thm}

\begin{Prf}
The second assertion follows from the first and from Prop.\ \ref{I'}. If
$I$ is finite, the \tg y induced, being $T_1$, must be discrete. In
particular for any $i\in I$\,\,$I\setminus\{i\}$ is closed hence $E$ contains
a \fn\ vanishing outside $i$ but not in $i$, thus contains the standard basis
and therefore is the whole $K^I$.   
\end{Prf}

\section{The spaces of operating \fn s, operating semigroups}

Let $E$ be a unital \fn\ space stable w.r.t.\ some non-additively-affine operation
$f:K^{n_0}\to K$. For $n=1,2,\ldots$ we have the spaces $G_n\subset K^{K^n}$
of the \fn s from $K^n$ to $K$ which operate on $E$. Clearly, any member
of $G_m$ operates on $G_n$ for any $m,n$, so they induce topologies
on $K^n$, $n=1,2,\ldots$, for which all members of the $G_m$ or $E$ are
\cs\ and also all ``vector-\fn s'' whose ``components'' are in $G_m$ or $E$
are \cs. Note that the topology on $K^{m+n}$ is {\em finer} than the product
topology $K^m\times K^n$, so the continuity of vector-\fn s is {\em not} a
consequence of the continuity of each component.
\parsk

These $G_n$ are an instance of what we call an {\EM operating semigroup}
defined as a sequence of linear subspaces $G_n\subset K^{K^n}$,
$n=1,2,\ldots$ such that $G_n$ contains the constants and the $n$ coordinate
\fn s $K^n\to K$ and every member of $G_m$ operates on $G_n$ for each $m,n$.
The minimal operating semigroup is the sequence of the spaces of affine
(i.e.\ constants $+$ linear \fn als). A possibly bigger one consists of all
additively affine \fn s. If some $G_n$ contains a non-additively-affine \fn,
the operating semigroup will induce \tg ies on the $K^n$ as above.
\parsk

If $G_n$ and $G_n'$ are operating semigroups, so is the sequence
$G_n\cap G_n'$. Therefore if $F$ is any set of \fn s $f:K^n\to K$
(possibly with different $n$ for different $f$'s), then there is a minimal
operating semigroup containing $F$ which
we call {\EM the operating semigroup generated by $F$} and denote by $G(F)$.
Since the \fn s operating on a \fn\ space form an operating semigroup,
any \fn\ space stable under $F$ is stable under $G(F)$.
\parsk

\section{Examples}

\begin{Prop}
Let $G_n$ be an operating semigroup $G_n$ which contains a non-additively-affine \fn.
Then (i)$\Rightarrow$(ii)$\Rightarrow$(iii):
\begin{itemize}
\item[(i)]
For all $n$, the \tg y induced on $K^n$ is the product
\tg y of the one induced on $K$.
\item[(ii)]
The \tg y induced on $K^2$ is the product \tg y of the one
induced on $K$.
\item[(iii)] 
The \tg y induced on $K$ is Hausdorff.
\end{itemize}
\end{Prop}

\begin{Prf}
(i)$\Rightarrow$(ii) is trivial, and (ii)$\Rightarrow$(iii) follows from the
fact that the diagonal in $K^2$ is always closed, being defined by a linear
equation. 
\end{Prf}

Note that in general (iii) does not imply (i). As an example take
$K=\bR(u)$ -- the field of real rational \fn s on one variable, and
$G_n$ -- the set of \fn s that are \cs\ on every \bd\
$\bR$-finite-dimensional set into a \bd\ $\bR$-finite-dimensional set.
This is an operating semigroup. The \tg y it induces on $K^n$ is the
finest \tg y that gives the usual \tg y on any $\bR$-finite-dimensional
subspace. It is Hausdorff, and as is well-known, (i) is violated. 
\parsk

\begin{Prop}
Let $G_n$ be an operating semigroup $G_n$ which contains a non-additively-affine \fn.
Then (i)$\Rightarrow$(ii):
\begin{itemize}
\item[(i)] 
The \tg y induced on $K$ is not Hausdorff.
\item[(ii)]
For all $n$, any two non-empty open sets, w.r.t.\ the \tg y induced on $K^n$, intersect.
\end{itemize}
\end{Prop}

\begin{Prf}
Firstly, if $n>1$ and $\exists$ two non-intersecting non-empty open sets
$U,V$ in $K^n$ then $\exists$ two such sets in $K^{(n-1)}$. Indeed, suppose not.
Then any $K^{(n-1)}$-section  of $K^n$ does not intersect either $U$ or $V$, which means that
projections of $U$ and $V$ on the first coordinate are disjoint, thus, being open
(as union of all parallel-to-1st-coordinate sections -- we don't need product \tg y),
one of them is empty hence either $U$ or $V$ is empty.
\parsk

Therefore it suffices to prove that any two non-empty open sets in $K$ intersect.
But if not, then some $x\ne y$ in $K$ would have disjoint \nd s, and since
the topology is invariant w.r.t.\ all invertible affine transformations of $K$, every
two points would have disjoint \nd s, i.e.\ the \tg y is Hausdorff.
\end{Prf}

The case $K=$ $\bR$ (the real numbers) or $\bC$ (the
complex numbers) is studied, in the context of \cs\ operations, in
\cite{deLeeuwKatznelson}, \cite{Sternfeld_A},
\cite{Sternfeld_DS} and \cite{SternfeldWeit}.
As is shown there, for $K=\bC$, if $E$ consists of \bd\ \fn s and is closed
w.r.t.\ the $\sup$-norm, then if some non-affine \cs\ \fn\ operates on $E$
then every analytic \fn\ operates on $E$: $E$ is an algebra, and if further
some non-analytic \cs\ \fn\ (say, complex conjugation) operates on $E$ then
any \cs\ \fn\ does so: $E$ is a self-adjoint algebra.%
\footnote{This is proved by showing that if the \fn\ space $E$ is stable
w.r.t.\ some \cs\ \fn\ $h$, then, being stable also w.r.t.\ the shifts of $h$
hence w.r.t.\ its convolution with, say, smooth \fn s (so one may assume $h$
smooth), one can make a derivative-like limiting process which shows that if
$h$ is not affine then $E$ is stable w.r.t.\ the square \fn, hence is an
algebra, and if $h$ is not analytic -- w.r.t.\ complex conjugation, hence is
self-adjoint.}
Thus in the latter case $G_m$, $m=1,2,\ldots$, the spaces of \fn s that
operate on $E$, contain all \cs\ \fn s, and the topologies that they induce
on $\bC^m$ are finer than the usual complex topologies. 
\parsk

They can be strictly finer. Indeed, if $I$ is infinite and $E$ is the space
of all $\bC$-valued \bd\ \fn s, then $G_m$ is the set of all
$f:\bC^m\to\bC$ that map every \bd\ sequence to a \bd\ sequence,
equivalently are \bd\ on every \bd\ set. The topology induced on
$\bC^n$ is the discrete one. If instead $E$ is the space of all $\bC$-valued
\fn s on $I$ then $G_m$ is the space of all \fn s $\bC^m\to\bC$, while the
topology induced on $\bC^m$ is again the discrete one.
\parsk

\begin{Rmk} \label{Rmk:lttc}
To quote another example, let $K$ be any ordered field, and suppose $E$ is
a lattice with pointwise lattice operations, which can be expressed as:
the \fn\ $a\mapsto a_+:=\max(a,0)$ from $K\to K$ operates on $E$. In this
case the \tg y induced on $K$ is finer than the order \tg y, hence
is Hausdorff, and of course the \tg ies on $K^n$ are finer than the
products of the order-\tg y.
\end{Rmk}

\section{Homomorphisms w.r.t.\ non-affine operations}

Fix a set $F$ of \fn s $f:K^n\to K$ (each possibly with different $n$), such
that $F$ contains a non-additively-affine \fn. We have the category of unital \fn\ spaces
$E\subset K^I$ stable under $F$ (equivalently, stable under $G(F)$, hence we
may assume at the outset that $F$ is an operating semigroup as defined above),
with {\EM unital $F$-homomorphisms}, i.e.\ $K$-linear mappings $\varphi$ which
commute with the action of every $f\in F$ and preserve $\bf1$
(thus preserve all constants).
\parsk

A unital $F$-homomorphism from a $E\subset K^I$ to $E'\subset K^{I'}$, both
unital and stable under $F$, has a graph, which may be viewed as a unital
\fn\ space $\varphi\subset K^{I\coprod I'}$, stable under $F$.
The \fn\ spaces $E$, $E'$ and $\varphi$ induce topologies
on $I$, $I'$ and $I\coprod I'$, resp. $E$ coincides with the space of
restrictions to $I$ of members of $\varphi$, therefore the \tg y on $I$ is
that of a subspace of $I\coprod I'$. Moreover, the only member of $\varphi$
vanishing on $I$ is identically $0$. Consequently $I$ is dense in
$I\coprod I'$. The \tg y of $I'$ as a subspace of $I\coprod I'$ is that
induced by the image of $\varphi$ and is coarser than that induced by $E'$.
\parsk

Thus we have

\begin{Cor}
Let $\varphi$ be a unital $F$-homomorphism from an $E\subset K^I$ to an
$E'\subset K^{I'}$, both unital and stable under $F$, and suppose $F$
contains a non-additively-affine member. Then for any $n=1,2,\ldots$ and $S\st K^n$
closed w.r.t.\ the \tg y induced by $G(F)_n$, if $\xi_1,\ldots\xi_n\in E$
assume together values in $S$ on $I$ then so do
$\varphi(\xi_1),\ldots,\varphi(\xi_n)$ on $I'$.
\end{Cor}

(Of course no such thing holds for general linear mappings between \fn\ spaces!)
\parsk

For instance, by Remark \ref{Rmk:lttc}, we obtain a neat proof of the
following fact:

\begin{Cor}
For any ordered field $K$, and for any lattice-homomorphism $\varphi$,
preserving the constants, between two $K$-valued \fn\ spaces $E$ and $E'$
which are lattices w.r.t.\ pointwise lattice operations, the image in $K^n$
of $(\varphi(\xi_1),\ldots,\varphi(\xi_n))$, $\xi_1,\ldots\xi_n\in E$ is
contained in the closure (w.r.t.\ the order-\tg y) of the image of 
$(\xi_1,\ldots,\xi_n)$.
\end{Cor}

When the \tg y induced by $G(F)$ is Hausdorff we have:

\begin{Cor}
Suppose $F$ is a set of \fn s $K^n\to K$ (possibly with different $n$ for
different \fn s), such that $F$ contains a non-additively-affine \fn\ and the \tg y
that $G(F)$ induces on $K$ is {\em Hausdorff}. Then any unital $F$-homomorphism
can be described as follows: take a unital $F$-stable \fn\ space $E\subset K^I$.
It induces a topology on $I$. Embed $I$ as a dense subset of some \tg ical
space $I\coprod I'$, such that each $\xi\in E$ has an extension (necessarily
unique!) to a \fn\ $\tilde\xi:I\coprod I'\to K$ \cs\ when $K$ is endowed with
the \tg y induced by $G(F)$, and map $\xi$ to $\tilde\xi|_{I'}$.
\end{Cor}

\end{document}